\documentclass[12pt]{amsart}
\usepackage{amscd,amssymb}
\usepackage[graph,frame,poly,arc]{xy}  
\usepackage[plainpages,backref,urlcolor=blue]{hyperref}
\usepackage{mathtools}

\theoremstyle{plain}
\newtheorem{thm}[subsection]{Theorem}

\newtheorem{prop}[subsection]{Proposition}
\newtheorem{cor}[subsection]{Corollary}

\theoremstyle{definition}
\newtheorem{rk}[subsection]{Remark}
\newtheorem{definition}[subsection]{Definition}
\newtheorem{ex}[subsection]{Example}

\newtheorem{question}[subsection]{Question}

\numberwithin{equation}{section}
\newcommand{\OO}{{\mathcal O}}

\newcommand{\Z}{\mathbb{Z}}

\newcommand{\C}{\mathbb{C}}
\newcommand{\PP}{\mathbb{P}}

\newcommand{\p}{\mathfrak {p}}

\DeclareMathOperator{\rank}{rank}
\DeclareMathOperator{\im}{im}
\DeclareMathOperator{\coker}{coker}

\DeclareMathOperator{\codim}{codim}

\DeclareMathOperator{\depth}{depth }

\DeclareMathOperator{\hit}{height}
\DeclareMathOperator{\gr}{grade}

\begin{document}

\title [Bourbaki modules and the module of Jacobian derivations]
{Bourbaki modules and the module of Jacobian derivations of projective hypersurfaces}

\author[Alexandru Dimca]{Alexandru Dimca}
\address{Universit\'e C\^ ote d'Azur, CNRS, LJAD, France and Simion Stoilow Institute of Mathematics,
P.O. Box 1-764, RO-014700 Bucharest, Romania}
\email{Alexandru.Dimca@univ-cotedazur.fr}

\author[Gabriel Sticlaru]{Gabriel Sticlaru}
\address{Faculty of Mathematics and Informatics,
Ovidius University
Bd. Mamaia 124, 900527 Constanta,
Romania}
\email{gabriel.sticlaru@gmail.com }



\begin{abstract} 
Two properties of projective hypersurfaces related to the module of Jacobian derivations, namely being tame and being plus-one generated,
are discussed in this paper. Tame hypersurfaces are related to Bourbaki ideals, and free hypersurfaces are the simplest examples of tame hypersurfaces. Plus-one generated hypersurfaces are the non free hypersurfaces which are closest to the free ones, and it is an open question whether all of them are tame.
\end{abstract}
 
\maketitle


\section{Introduction} 

Let $S=\C[x_0,x_1, \ldots, x_n]$ be the polynomial ring in $n+1 \geq 3$ variables $x_0,x_1, \ldots, x_n$ with complex coefficients, and let $V:f=0$ be a reduced hypersurface of degree $d\geq 3$ in the complex projective plane $\PP^n$. 
We denote by $J_f$ the Jacobian ideal of $f$, i.e., the homogeneous ideal in $S$ spanned by the partial derivatives $f_j$ of $f$ with respect to $x_j$ for $j=0, \ldots,n$, and  by $M(f)=S/J_f$ the corresponding graded quotient ring, called the Jacobian (or Milnor) algebra of $f$.
Consider the graded $S$-module of Jacobian syzygies of $f$ or, equivalently, the module of derivations killing $f$, namely
\begin{equation}
\label{eqD0}
D_0(f)= \{\theta \in Der(S) \ : \ \theta(f)=0\}= \{\rho=(a_0,a_1, \ldots, a_n) \in S^{n+1} \ : \  \sum_0^na_jf_j=0\}.
\end{equation}
This module is also denoted in the literature by 
$AR(f)$ (i.e., all Jacobian relations for $f$) or by ${\rm Syz}(f)$ (i.e., the Jacobian syzygies of $f$). 
We say that $V:f=0$ is an {\it $m$-syzygy hypersurface} if  the module $D_0(f)$ is minimally generated by $m$ homogeneous syzygies, say $\rho_1,\rho_2,\ldots ,\rho_m$, of degrees $d_j=\deg \rho_j$ ordered such that $$d_1\leq d_2 \leq \ldots \leq d_m.$$ 
In this note we assume  that $V$ is not the cone over a hypersurface in $\PP^{n-1}$, hence $d_1 \geq 1$.
We call these degrees $(d_1, \ldots, d_m)$ the {\it exponents} of the hypersurface $V$ and $\rho_1,\ldots,\rho_m$ a {\it minimal set of generators } for the module  $D_0(f)$. 
The smallest degree $d_1$ is sometimes denoted by ${\rm mdr}(f)$ and is called the minimal degree of a Jacobian relation for $f$. 

Since the $S$-module $D_0(f)$ has rank $n$, it follows that one always has
$$m \geq n.$$
The hypersurface $V$ is {\it free} when $m=n$, since then  $D_0(f)$ is a free module of rank n, see for instance \cite{freeS, KS,Sim2,ST,To}. Let $\Sigma$ be the singular subscheme of $V$, defined by the Jacobian ideal. Then $s=\dim \Sigma \leq n-2$ since $V$ is reduced, and note that the Hilbert polynomial 
$H(M(f),t)$ of the graded $S$-module $M(f)$ has the form
\begin{equation}
\label{eqD1}
H(M(f),t)= \frac{\deg(\Sigma)}{s!}t^s+ \ldots
\end{equation}
where $\ldots $ stand for lower degree terms in $t$ and $\deg(\Sigma)$ is by definition the degree of the scheme $\Sigma$, or of the ideal $J_f$.

Beyond the free hypersurfaces, there is the class of {\it plus-one generated } hypersurfaces, for short POG hypersurfaces, defined in the case of hyperplane arrangements by Takuro Abe in \cite{Abe}.
With our notation, {\it choosing  an order on the syzygies $\rho_j$'s, but without assuming that the exponents $d_j$'s are ordered}, a POG hypersurface is a hypersurface $V:f=0$ in $\PP^n$ such that $m=n+1$ and the second order syzygies of $f$, that is the relations among the minimal generators $\rho_j$ are generated by just one relation, namely
\begin{equation}
\label{eqD2}
b_1\rho_1+ \ldots + b_n \rho_n + \ell \rho_{n+1}=0,
\end{equation}
where $\ell \in S_1$ is a linear form. 
If we want to emphasize that this definition depends on the choice of the generators, we say that $V:f=0$ is a POG hypersurface with respect to 
$$\overline \rho=(\rho_1, \ldots, \rho_{n+1}).$$ Since the set of generators is minimal, it follows that none of the coefficients $b_j$ is a non-zero constant, that is either $b_j=0$ or $b_j \ne 0$ and $\deg b_j \geq 1$.
The hypersurface $V$ is called {\it strictly plus-one generated } if in addition $\ell \ne 0$. Note that in this case, we get $d_{n+1} \geq d_j$ for all $j$ such that $b_j \ne 0$, even if we do not order the degrees $d_j$.
When $n=2$, any POG curve is in fact a strictly POG curve, see \cite[Corollary 2.2]{minTjurina} for a more general result applying to any 3-syzygy curve, but in higher dimensions this nice property may fail, see Example \ref{ex10} and Question \ref{q3}.

A POG hypersurface with respect to 
$\overline \rho$ such that $d_1 \leq \ldots \leq d_{n}=d_{n+1}$ is called a {\it nearly-free hypersurface}, see \cite[Definition 5.3]{freeS} in the case $n=3$. 

In this note we first collect in Section 2 some basic facts on the module $D_0(f)$ in the case of curves in $\PP^2$, that is $n=2$. In particular, we list some properties for
free and plus-one generated curves. We will try to extend these properties in the fourth section to the case $n>2$.

 The fact that the case $n>2$ is much harder to understand comes perhaps from the following technical point.
Assume we fix 
$$ \hat \rho=(\rho_1, \rho_2, \ldots, \rho_{n-1})$$
a subset of $(n-1)$ elements in a minimal set of generating syzygies for $f$, without an ordering of the degrees. Then we have a sequence of graded $S$-modules
\begin{equation}
\label{eqB1}
 0 \to S(-d_1) \oplus \ldots \oplus S(-d_{n-1})  \xrightarrow{u} D_0(f) \xrightarrow{v} B(f,\hat \rho)(d_1+\ldots +d_{n-1}-d+1) \to 0,
\end{equation}
where 
$$u(a_1,\ldots, a_{n-1})= \sum_{j=1}^{n-1}a_j\rho_j$$
and the morphism $v$ is defined as follows. For any syzygy $\rho \in D_0(f)$, consider the $(n+1) \times (n+1)$-matrix $M(\hat \rho, \rho)$ where the first row is
given by
$$(x_0, x_1, \ldots, x_n),$$
which are the components of the Euler derivation $E$,
the $j$-th row, for $2 \leq j \leq n$ is given by the components of the $j-1$-syzygy $\rho_{j-1}$ fixed above, and the $n+1$-st row by the components of $\rho.$ Then it is easy to show that the polynomial $\det M(\hat \rho, \rho)$ is divisible by $f$ and we set
$$v(\rho)= \frac{\det M(\hat \rho, \rho)}{f}.$$
Then the image of $v$ is a (shifted) ideal in $S(d_1+\ldots +d_{n-1}-d+1)$,
which we denote by $B(f,\hat \rho)=\im v(-(d_1+\ldots +d_{n-1}-d+1))$, and which is called the Bourbaki ideal of the module $D_0(f)$, see \cite{SS,DStJump,JNS} for the case $n=2$ and \cite{DmaxS} for the case $n=3$.
 Following the terminology introduced in \cite{DmaxS}, we give the following.
\begin{definition}
\label{defT}
We say that the hypersurface $V:f=0$ is {\it tame with respect to the subset $\hat \rho$} if the sequence \eqref{eqB1} is exact.
\end{definition}
\begin{rk}
\label{rkT}
In the literature, there is a notion of {\it tame arrangements or tame hypersurfaces}, see \cite{TY} or \cite{Wa}, and another notion of {\it tame polynomials}, see \cite{Bt} or \cite{NZ}. {\it Our notion of tameness does not seem to be related to either of these notions.}
Note also that only for hypersurfaces $V:f=0$ which are tame with respect to the subset $\hat \rho$, the ideal $B(f,\hat \rho)$ is the Bourbaki ideal of the $S$-module $D_0(f)$ in the sense of the {\it general theory of Bourbaki ideals}, see \cite{HKS} and our discussion in Section 3 below.
\end{rk}

In Section 3 we compare our sequence \eqref{eqB1} to the usual notion of a Bourbaki sequence for a (graded) $S$-module. Theorem \eqref{thmT0} is the main result here, giving necessary and sufficient conditions for the sequence \eqref{eqB1} to be exact. The proof of this result uses a key result due to J. Herzog, S. Kumashiro and D. Stamate
in \cite{HKS}.
As we shall see, being tame with respect to $\hat \rho$  depends on the choice of $\hat \rho$, see
Example \ref{ex01}. 

In Section 4 we show that 
for tame hypersurfaces, many properties of the $S$-module $D_0(f)$ which hold for $n=2$ extend to the case $n>3$, see Theorems \ref{thmFF} and \ref{thmF0}, but fail as soon as we omit this property, see Remark \ref{rkFF} and Remark \ref{rkF0}.

We end the paper with a number of natural questions concerning 
the module of Jacobian derivations, which might foster new research in this area.

\section{Basic facts  to Jacobian syzygies of curves } 

When $n=2$, we use the simpler notations $x_0=x, \ x_1=y$ and $x_2=z$.
According to Hilbert Syzygy Theorem, the general form of the minimal resolution for the Milnor algebra $M(f)$ of a curve $C:f=0$ that is assumed to be not free, is the following
\begin{equation}
\label{res2A}
0 \to \oplus_{i=1} ^{m-2}S(-e_i) \to \oplus_{j=1} ^mS(1-d-d_j)\to S^3(1-d)  \to S,
\end{equation}
with $e_1\leq \ldots \leq e_{m-2}$ and $d_1\leq \ldots \leq d_m$.
The last morphism in \eqref{res2A} is given by
$$(a,b,c) \mapsto af_x+bf_y+cf_z,$$
and has as kernel precisely the graded $S$-module of Jacobian syzygies $D_0(f)$ shifted by $(1-d)$.
It follows from \cite[Lemma 1.1]{HS} that one has
\begin{equation}
\label{res2B}
e_j=d+d_{j+2}-1+\epsilon_j,
\end{equation}
for $j=1,\ldots,m-2$ and some integers $\epsilon_j\geq 1$. With this notation, the resolution in \eqref{res2A} yields the following minimal resolution for the graded $S$-module $D_0(f)$:
\begin{equation}
\label{res2A1}
0 \to F_1=\oplus_{i=1} ^{m-2}S(-e'_i) \to F_0=\oplus_{j=1} ^mS(-d_j),
\end{equation}
where $e'_j=d_{j+2}+\epsilon_j$, for $j=1,\ldots,m-2$. Therefore, the integers $\epsilon _j$ gives the shifts of the degrees of the generators when we pass from $F_0$ to $F_1$ in the minimal resolution of $D_0(f)$.
Using \cite[Formula (13)]{HS}, it follows that one has
\begin{equation}
\label{res2C}
d_1+d_2=d-1+\sum_{i=1} ^{m-2}\epsilon_j \geq d-1+m-2,
\end{equation}
which gives in particular the following bounds on the number $m$ of a minimal generating set
\begin{equation}
\label{inq1}
2 \leq m \leq d_1+d_2-d+3,
\end{equation}
see \cite[Proposition 2.1]{maxTjurina} and \cite{E2}.
On the other hand,
for a reduced degree $d$ curve $C$, one has 
\begin{equation}
\label{e1}
m \leq d+1,
\end{equation}
see \cite[Proposition 2.1]{maxTjurina}, while for the case $C$ a line arrangement, we have the slightly stronger inequality
\begin{equation}
\label{e2}
 m \leq d-1,
 \end{equation}
 see \cite[Corollary 1.3]{DStSat},
It is known that, for a reduced degree $d$ curve $C$, one has 
\begin{equation}
\label{e3}
d_m \leq 2d-4,
\end{equation}
 see \cite[Corollary 11]{CD}, and for the case $C$ a line arrangement, one has the much stronger inequality
 \begin{equation}
\label{e4}
 d_m \leq d-2,
 \end{equation}
 see \cite[Corollary 3.5]{Sch0}.

Next we collect some characterizations of free curves. The first one is very easy to check, as soon as we have the first two exponents of $C$.
\medskip

\noindent  (C1) \  \   A free curve $C$ is  characterized by the condition $d_1+d_2=d-1$, 
 
\medskip

\noindent and hence in this case
$r=d_1 <d/2$.
For a proof,  see \cite{ST}. 
Recall that the total Tjurina number $\tau(C)$ of a reduced plane curve $C:f=0$ is just the degree $\deg \Sigma$ of the singular subscheme, or equivalently the sum of the Tjurina numbers of all the singularities of $C$.
The main result in \cite{duPCTC}, restated below in
Theorem \ref{thmCTC}, gives an upper  bound for the total Tjurina number $\tau(C)$ of $C$ as a function of $d$ and $r=d_1$. Moreover, we know that when $r<d/2$ the upper bound is obtained exactly for a free curve, see Corollary \ref{corCTC} below. Moreover, it was shown in \cite{expo} that for any $r<d/2$ there is a free curve with exponents $d_1=r$ and $d_2=d-1-r$. 

\begin{thm}
\label{thmCTC}
For positive integers $d$ and $r$, set $
\tau(d,r)_{max}= (d-1)^2-r(d-r-1).$
Then, if $C:f=0$ is a reduced curve of degree $d$ in $\PP^2$ and  $r=d_1={\rm mdr}(f)$,  one has
$$\tau(C) \leq \tau(d,r)_{max}.$$

\end{thm}
For a proof, see \cite[Theorem 3.2]{duPCTC} as well as \cite{E} for an alternative approach.

At the end of the proof of Theorem \ref{thmCTC}, in \cite{duPCTC}, the authors state the following very interesting consequence (of the proof, not of the statement) of Theorem \ref{thmCTC}.
\begin{cor}
\label{corCTC} Let $C:f=0$ be a reduced curve of degree $d$ in $\PP^2$ and  $r=d_1={\rm mdr}(f)$. One has
$$ \tau(C) =\tau(d,r)_{max}$$
if and only if $C:f=0$ is a free curve, and then $r <d/2$.
\end{cor}
In the paper \cite{Dmax}, there is an alternative proof of Corollary \ref{corCTC}. Moreover, a plane curve $C$ is nearly free, which can be defined by the property $m=3$, $d_1+d_2=d$ and $d_2=d_3$ if and only if a similar property holds. Namely, one has the following result, see \cite{Dmax}.

\begin{cor}
\label{corCTC2} Let $C:f=0$ be a reduced curve of degree $d$ in $\PP^2$ and  $r={\rm mdr}(f)$.
One has
$$ \tau(C) =\tau(d,r)_{max}-1$$
if and only if $C:f=0$ is a nearly free curve, and then $r  \leq d/2$.
\end{cor}

Next we state a characterization of POG curves, which is very easy to check, as soon as we have the first two exponents of $C$.
\medskip

\noindent   (C2) \  \   A non free curve $C$  satisfies $d_1+d_2 \geq d$. Moreover, the equality holds if and only if $C$ is a POG curve,
 
\medskip
see \cite[Theorem 2.3]{minTjurina}.

The exact sequence \eqref{eqB1} in the case $n=2$ goes back to
\cite{duPCTC}, and was gradually developed in \cite{Dmax,SS,DStJump,JNS}. The main fact is the following result.

\medskip

\noindent   (C3) \  \   Any curve $C:f=0$ in $\PP^2$ is tame with respect to  $\hat \rho=(\rho_1)$, a syzygy of minimal degree in $D_0(f)$,

 \medskip
 
see \cite[Lemma 3.2 and the exact sequence (3.3)]{Dmax}.
\section{Tame hypersurfaces vs. Bourbaki exact sequences} 
The tameness property plays a key role in this section, so we start with
a discussion of this property. First we recall the notion of Bourbaki sequence for a module.
Let $R$ be a Noetherian ring and $M$ an $R$-module, where both $R$ and $M$ are not supposed to be graded. Then a Bourbaki sequence of $M$ means a short exact sequence
\begin{equation}
\label{eqBT}
0 \to F \to M \to I \to 0
\end{equation}
of $R$-modules, where $F$ is a free $R$-module and $I$ is an ideal of $R$. It is known that a Bourbaki sequence of $M$ always exists if $R$ is a normal domain and $M$ is a finitely generated torsionfree $R$-module, see \cite[Chapter VII, §4, 9. Theorem 6]{B}. Bourbaki sequences were considered also for graded modules. For instance, \cite[Theorem 2.1 and Lemma 2.2 ]{HKS} yields the following result.
\begin{thm}
\label{thmT1}
Let $S$ be the polynomial  ring $\C[x_0, \ldots, x_n]$ with $n \geq 2$ as above and $M$ a finitely generated reflexive graded $S$-module of rank $r>1$. Then there is a graded Bourbaki sequence 
$$0 \to S(-D)^{r-1} \to M \to I (s)\to 0,$$
where $D$ is the maximal degree of a generator of $M$ in a minimal set of generators for $M$, $s$ is an integer and $I$ a graded ideal in $S$ with $\gr (I) =2$.
\end{thm}
Note that the sequence \eqref{eqB1} is not a Bourbaki sequence for $D_0(f)$ as in Theorem \ref{thmT1}, not even when $n=\rank D_0(f)=2$, since in this case we usually choose  $\hat \rho=(\rho_1)$,
with $\rho_1$ a minimal degree generator for $D_0(f)$. When $n>2$,
our sequence  \eqref{eqB1}  relates in a more precise way the numerical invariants of the hypersurface $V:f=0$, assumed to be tame with respect to $\hat \rho$, to the degrees of the of the syzygies in 
$\hat \rho$, see for instance Theorems \ref{thmFF}, \ref{thmF0} and \ref{thmPOG20} below. However, Theorem \ref{thmT1} clearly implies the following result.
\begin{cor}
\label{corT1}
Any reduced hypersurface $V:f=0$ in $\PP^n$ is tame for some subset
$\hat \rho$ if the minimal set of generators for $D_0(f)$ satisfies
$$d_1=d_2= \ldots = d_m.$$
\end{cor}
In general, the sequence \eqref{eqB1} is exact at the first term, that is $u$ is injective exactly when the $(n-1) \times (n+1)$-matrix $M(\hat \rho)'$ which has as rows the components of the syzygies in $\hat \rho$ satisfies
\begin{equation}
\label{eqB2}
\rank M(\hat \rho)'=n-1,
\end{equation}
where the rank is computed say over the field of fractions $K$ of $S$.
To check whether the sequence \eqref{eqB1} is exact at the second term, that is 
$$\im u= \ker v,$$
it is more subtle. We have the following result.
\begin{thm}
\label{thmT0}
Let $I(\hat \rho)$ be the ideal in $S$ generated by all the $(n-1)$-minors in the matrix $M(\hat \rho)'$. Then the sequence \eqref{eqB1} is exact
if and only if
$$\hit (I(\hat \rho) \geq 2.$$
\end{thm}
\proof
The proof uses the key result which appear in \cite[Theorem 1.2]{HKS}, and covers the ungraded situation described in \eqref{eqBT}.
To apply this result, we recall the following facts. First the polynomial ring $S$ is clearly a normal domain of dimension $n+1 \geq 2$, since we have assumed in this paper that $n \geq 2$. The $S$-module $D_0f)$ is finitely generated and torsion free, and has rank $n \geq 2$.
Moreover, $D_0(f)$ is reflexive and sits in the exact sequence of $S$-modules
$$0 \to D_0(f) \to S^{n+1} \to J_f(d-1) \to 0,$$
which satisfies the hypotheses of \cite[Theorem 1.2 (a)]{HKS}. Then
 \cite[Theorem 1.2 (a)]{HKS} tells us that $\hit (I(\hat \rho) \geq 2$ if and only if the following sequence of (ungraded) $S$-modules
\begin{equation}
\label{eqB20} 
  0 \to S(-d_1) \oplus \ldots \oplus S(-d_{n-1})  \xrightarrow{u} D_0(f) \to \coker (u) \to 0
 \end{equation} 
  is exact and $I'=\coker (u)=D_0(f)/\im(u)$ is an ideal in $S$. Since the grading plays no role here, we may replace
 $S(-d_1) \oplus \ldots \oplus S(-d_{n-1}) $ by $S^{n-1}$ to have the formulation exactly as in  \cite[Theorem 1.2 (a)]{HKS}.
 In fact, the injectivity of $u$ is clear, since the condition \eqref{eqB2} fails if and only if $I(\hat \rho)=0$, which contradicts the inequality
$\hit (I(\hat \rho) \geq 2$, since $\hit (0)=0$.
If the equality  $\im u= \ker v$ fails, then any element $\rho \in \ker v$
such that $\rho \notin \im v$ gives rise to a torsion element in $I'$, which is a contradiction with the fact that $I'$ is an ideal in $S$. This implies  by \cite[Theorem 1.2 (a)]{HKS} that in this case $\hit (I(\hat \rho)) <2$. If the equality $\im u= \ker v$ holds,
then the exact sequence \eqref{eqB1} may be identified to the sequence
\eqref{eqB20}, and we conclude, using again \cite[Theorem 1.2 (a)]{HKS}, that in this case $\hit (I(\hat \rho) \geq 2$.
This completes our proof, 
since the sequence \eqref{eqB1} is exact at the third term by the definition of the Bourbaki ideal. 

\endproof

\begin{rk}
\label{rkT0}
Recall that for any ideal $I \subset S$, one has
$$\hit(I)= \min_{I \subset \p}\hit(\p)=\min_{I \subset \p}\codim V(\p),$$
see for instance \cite[Theorem A.16]{BH},
where the minimum is over all prime ideal $\p$ containing $I$, in other words over all the irreducible varieties $V(\p)$ contained in the scheme $V(I)$ defined by $I$ in $\C^{n+1}$. On the other hand, we clearly have
$$\codim V(I)=\min_{I \subset \p}\codim V(\p).$$
Hence the condition
$\hit(I(\hat \rho))\geq 2$ in Theorem \ref{thmT0} is equivalent to
$$\codim V(I(\hat \rho)) \geq 2.$$
In other words, $V(I)$ contains no codimension 1 component, that is the set of all the $(n-1)$-minors in the matrix $M(\hat \rho)'$ are not divisible by any polynomial in $S$ of degree $>0$.
Similarly, since in the polynomial ring $S$ we have
$$\gr (I)= \depth (I)= \hit (I)$$
for any ideal $I$, it follows that the condition $\gr (I)=2$ in Theorem \ref{thmT1} is equivalent to 
$$\codim I=2.$$
\end{rk}

Now it turns out that for $n=2$, any reduced plane curve is tame, see the property (C3) in Section 2. For $n=3$, the first condition \eqref{eqB2} holds for any surface and any choice $\hat \rho=(\rho_1, \rho_2)$ with $\rho_1$ of minimal degree and $\rho_2 \notin S\cdot \rho_1$. For $n \geq 4$, the tameness condition is even more subtle.
Here is an example illustrating this property in the case $n=3$, see also Remark \ref{rkF0} and Example \ref{ex1} below.  We use in this paper the simpler notation $x_0=x, \ x_1=y, \ x_2=z$ and $x_3=w$ when $n=3$.

\begin{ex}
\label{ex01} The following example is related in particular to Corollary \ref{corT1}.
 Consider the surface
$$V:f=(x^2-zw)(x^2+y^2+z^2)=0,$$
which is a union of two quadratic cones in $\PP^3$.
Then the computer algebra software SINGULAR \cite{Sing} gives $m=7$ and the
exponents $(2,2,2,2,2,2,2)$. More precisely, we have
$$\rho_1=(0,-z^2,yz,-yw), \ \rho_2=(z^2,0,-xz,2xz+xw), \ \rho_3=
 (yw,-2xz-xw, 2xy,0), $$
$$\rho_4=(yz,-xz,0,2xy), \rho_5= (-xy,2x^2+y^2+z^2,0,-2yw),$$
$$\rho_6=(-2xz+xw,2yz+yw,4x^2+2z^2+zw,-6zw-3w^2)$$
 and $$
\rho_7=(4xz-xw,-yw,-4x^2-zw,8x^2+3w^2).$$
It follows that $z \rho_4=x\rho_1+y\rho_2$. Hence, if we choose the subset
$\hat \rho=(\rho_1,\rho_2)$, then $V$ is not tame with respect to $\hat \rho$, since $v(\rho_4)=0$ and clearly $\rho_4 \notin \im u$.
However, if we take the subset $\hat \rho=(\rho_1, \rho_6)$, a direct computation with SINGULAR shows that $\hit(I(\hat \rho))=2$, and hence $V$ is tame with respect to this $\hat \rho$ by Theorem \ref{thmT0}.
\end{ex}

\section{Results related to Jacobian syzygies of hypersurfaces } 

It is clear that any free hypersurface is tame, for any choice of the first $n-1$ generators $\hat \rho$. The following result can be regarded as an extension of property (C1) from previous section to dimensions $n>2$.

\begin{thm}
\label{thmFF} 
Let $V:f=0$ be a tame hypersurface with respect to the subset $\hat \rho=(\rho_1, \ldots, \rho_{n-1})$ and having the (unordered) set of exponents  $(d_1, \ldots, d_m)$ for some $m \geq n$. Then the following hold.
\begin{enumerate}

\item The Bourbaki ideal
$B(f,\hat \rho) \subset S$ is saturated and $\dim S/B(f,\hat \rho) \leq n-1$. Moreover, the minimal number of generators of $B(f,\hat \rho) $ is
$m-n+1$.

\item The hypersurface $V:f=0$ is free if and only if $B(f,\hat \rho) = S$.

\item The hypersurface $V:f=0$ is an $(n+1)$-syzygy hypersurface if and only if the ideal $B(f,\hat \rho)$ is a codimension 2 complete intersection.

\item  For any $j \geq n$ one has
$$d_1+ \ldots +d_{n-1}+d_j \geq d-1,$$
and the equality holds for some $j \geq n$ if and only if the tame hypersurface $V$ is free.
\end{enumerate}

\end{thm}
\proof
To prove (1),
let $\delta=d_1+ \ldots +d_{n-1} -( d-1).$ Then the exact sequence
\eqref{eqB1} can be written in the form
\begin{equation}
\label{eqF10}
0 \to S(-d_1-\delta) \oplus \ldots \oplus S(-d_{n-1}-\delta)  \xrightarrow{u} D_0(f)(-\delta) \xrightarrow{v} B(f,\hat \rho) \to 0.
\end{equation}
By taking the associated coherent sheaves, we get the following exact sequence
$$ 0 \to \OO(-d_1-\delta) \oplus \ldots \oplus \OO(-d_{n-1}-\delta)  \xrightarrow{u} \tilde D_0(f)(-\delta) \xrightarrow{v} \tilde B(f,\hat \rho) \to 0.$$
Looking at the corresponding long exact sequence of cohomology groups, and using the fact that
$$H^0(\PP^n, \tilde B(f,\hat \rho))=B(f,\hat \rho)^s,$$
the saturation of the ideal $B(f,\hat \rho)$, and
$$H^1(\PP^n, \OO(u))=0$$
for any $n \geq 2$ and any integer $u$, we get the exact sequence
\begin{equation}
\label{eqF11}
 0 \to S(-d_1-\delta) \oplus \ldots \oplus S(-d_{n-1}-\delta)  \xrightarrow{u} D_0(f)(-\delta) \xrightarrow{v} B(f,\hat \rho)^s \to 0.
\end{equation}
By comparing \eqref{eqF10} and \eqref{eqF11}, we see that 
$B(f,\hat \rho)=B(f,\hat \rho)^s$, hence the Bourbaki ideal is saturated.
To prove that $\dim S/B(f,\hat \rho) \leq n-1$ one can use \cite[Lemma 2.2]{HKS}, since the $S$-module $D_0(f)$ is reflexive by definition. The claim about the minimal number of generators follows from the exact sequence \eqref{eqF10}, which shows that any system of generators for $B(f,\hat \rho)$ lifts to a set of syzygies which together with the subset $\hat \rho$ generate $D_0(f).$
This remark completes the proof of (1) and also proves (3). Indeed, an ideal $I \subset S$ generated by two homogeneous polynomial such that $\dim S/I \leq n-1$ is a complete intersection of dimension $n-1$.
To prove (2) and (4), note that
 if $d_1+ \ldots +d_{n-1}+d_j < d-1$, then $e=\delta+d_j <0$. By taking the
 $e$-homogeneous component of the sequence \eqref{eqF10}, we get
 $v(\rho_j)=0$ since $S_e=0$. By tameness, this implies that $\rho_j \in \im u$, a contradiction with the minimality of the generating set $\rho_k$'s.
 If the equality holds, then $e=0$ and $v(\rho_j)$ is a non-zero constant,
 which implies that $B(f,\hat \rho)=S$. A well-known result by K. Saito \cite{KS}, see for instance \cite[Theorem 8.1]{DHA} for a formulation in our setting, shows that this happens if and only if $V$ is free.
 \endproof

\begin{rk}
\label{rkFF}  
Without the tameness hypothesis, there are non-free hypersurfaces such that
$$d_1+ \ldots +d_n =d-1,$$
see Example \ref{ex1} below. 
\end{rk}
The following result can be regarded as an extension of Corollary \ref{corCTC} to dimensions $n>2$. In the case $n=3$ this result is essentially \cite[Theorem 1.4]{DmaxS}.
\begin{thm}
\label{thmF0} 
Let $V:f=0$ be a free hypersurface in $\PP^n$ with exponents
$(d_1, \ldots, d_n)$ and $d= \deg f$. Then the singular subscheme $\Sigma$ has dimension $n-2$ and its degree is given 
by
$$ \deg \Sigma= (d-1)^2-s_2,$$
where  $s_2=\sum_{1\leq i<j\leq n}d_id_j$.
Moreover, let  $V':f'=0$ in $\PP^n$ be a tame hypersurface with respect to the subset
$\hat \rho'=(\rho_1',\ldots, \rho_{n-1}')$, of degree $d'$  and with exponents  $(d_1', \ldots, d_m')$ for some $m \geq n$. If the singular subscheme $\Sigma'$ of $V'$ has dimension
$n-2$, then
the degree of the singular subscheme $\Sigma'$ satisfies the inequality
$$ \deg \Sigma ' \leq (d'-1)^2-\sum_{1\leq i<j\leq n}d_i'd_j'$$
where in the last sum $d_n'$ is replaced by 
$$d_n''=d'-1-\sum_{j=1}^{n-1}d_j'.$$
The equality holds here if and only if $V'$ is free and hence in particular $m=n$.
\end{thm}

\proof
Since $V$ is a free hypersurface in $\PP^n$ with exponents
$(d_1, \ldots, d_n)$, it follows that
\begin{equation}
\label{eqF1}
 D_0(f)=\oplus_{j=1}^nS(-d_j)
 \end{equation}
and hence, for any large integer $k$ we have
$$\dim D_0(f)_k=\sum_{j=1}^n\binom{k+n-d_j}{n}.$$
The obvious exact sequence
\begin{equation}
\label{eqF2}
 0 \to D_0(f) \to S^{n+1} \to S(d-1) \to M(f)(d-1) \to 0
 \end{equation}
 implies that
 $$\dim D_0(f)_k=(n+1)\binom{k+n}{n}-\binom{k+n+d-1}{n}+\dim M(f)_{k+d-1}.$$
 Therefore we get
\begin{equation}
\label{eqF3} 
 \dim M(f)_{k+d-1}=\sum_{j=1}^n\binom{k+n-d_j}{n}-(n+1)\binom{k+n}{n}+\binom{k+n+d-1}{n}.
\end{equation} 
 Define polynomials $A_j(t) \in \Z[t]$ for $j=1, \ldots,n$ such that
 $$n!\binom{k+n}{n}=k^n+A_1(n)k^{n-1}+A_2(n)k^{n-2}+\ldots,$$
 where $\ldots$ denotes terms in $k$ with exponent $<n-2$. With this notation, one has
$$n!\binom{k+n-d_j}{n}=$$
$$k^n+(A_1(n)-nd_j)k^{n-1}+(A_2(n)-(n-1)d_jA_1(n)+\binom{n}{2}d_j^2)k^{n-2}+\ldots$$
and
 $$n!\binom{k+n+d-1}{n}=$$
$$k^n+(A_1(n)-n+nd)k^{n-1}+(A_2(n)+(n-1)(d-1)A_1(n)+\binom{n}{2}(d-1)^2)k^{n-2}+\ldots .$$
Using the equality \eqref{eqF3} it follows that 
$$n! \dim M(f)_{k+d-1}=n(d-1-\sum_{j=1}^nd_j)k^{n-1}+$$
$$+\left((n-1)(d-1-\sum_{j=1}^nd_j)A_1(n))+\binom{n}{2}((d-1)^2+\sum_{j=1}^nd_j^2)\right)k^{n-2}+ \ldots .$$
 Comparing this formula with the equality \eqref{eqD1}, we see that the coefficient of $k^{n-1}$ has to vanish, and hence we get
 $$d-1=\sum_{j=1}^nd_j,$$
 as we have already seen in Theorem \ref{thmFF}.
It follows that the coefficient of $k^{n-2}$ is strictly positive and hence
$s=\dim \Sigma=n-2$ and
$$\deg \Sigma=\frac{(n-2)!}{n!}\binom{n}{2} ((d-1)^2+\sum_{j=1}^nd_j^2)=\frac{(\sum_{j=1}^nd_j)^2+\sum_{j=1}^nd_j^2}{2}=s_1^2-s_2,$$
where $s_1=\sum_{j=1}^nd_j$ and $s_2=\sum_{1\leq i<j\leq n}d_id_j$.
Finally, consider a tame hypersurface $V':f'=0$. Then the sequence of graded $S$-modules \eqref{eqB1} is exact. This implies that
$$\dim D_0(f)_k \leq \sum_{j=1}^{n-1}\binom{k+n-d_j'}{n} + \binom{k+n-d'+1+\sum_{j=1}^{n-1}d_j'}{n},$$
with equality if and only if the corresponding Bourbaki ideal $B(f', \hat \rho ')$ is the polynomial ring $S$ by Theorem \ref{thmFF}. Using as above the exact sequence
\eqref{eqF2} it follows that
$$ \dim M(f')_{k+d'-1}\leq \sum_{j=1}^{n-1}\binom{k+n-d_j'}{n} + \binom{k+n-d'+1+\sum_{j=1}^{n-1}d_j'}{n}-$$
$$
 -(n+1)\binom{k+n}{n}+\binom{k+n+d'-1}{n}.$$
 Using the above notation, we see that
 $$n!  \dim M(f')_{k+d'-1} \leq (d'-1)^2-\sum_{1\leq i<j\leq n}d_i'd_j'$$
 where $d_n'$ is replaced by $d'-1-\sum_{j=1}^{n-1}d_j'.$
 
\endproof

\begin{rk}
\label{rkF0}  
1) Without the tameness hypothesis with respect to $\hat \rho'$, there are non-free hypersurfaces $V'$ such that the inequality for $\deg \Sigma'$
is an equality,
see Example \ref{ex1} below. Moreover, the inequality for $\deg \Sigma '$ holds even when $d_n'' \leq 0$, see Examples \ref{ex01} and \ref{ex0}.

2) Note that there are tame hypersurfaces $V$ such that the corresponding singular subscheme satisfies $\dim \Sigma <n-2$. For instance, the surface in $\PP^3$ given by 
$$V:f=xyz-w^3=0$$
has isolated singularities. Them $m=5$ and the generating syzygies are
$$\rho_1=(0,y,-z,0), \rho_2=(x,0,-z,0), \ \rho_3=(3w^2,0,0,yz),$$
$$\rho_4=(0,3w^2,0,xz) \text{ and } \rho_5=(0,0,3w^2,xy).$$
Then it is easy to see that $V$ is tame with respect to $\hat \rho=(\rho_1,\rho_2)$. Indeed, the matrix $M(\hat \rho)$ has rank 2 as required in \eqref{eqB2}. To show that $\ker v=\im u$ also holds, we proceed as follows. We can take $\rho=A\rho_3+B\rho_4+C\rho_5$
such that $A,B,C \in S$ and $v(\rho)=0$. A direct computation shows that this means
$$Ayz+Bxz+Cxy=0.$$
This implies $A=xA_1$, $B=yB_1$ and $C=zC_1$ with $A_1,B_1,C_1 \in S$. It follows that $A_1+B_1+C_1=0$. Therefore one has
$$\rho=A_1(x\rho_3-z\rho_5)+B_1(y\rho_4-z\rho_5).$$
Since both $(x\rho_3-z\rho_5)$ and $(y\rho_4-z\rho_5)$ are clearly in
$\im u$, it follows that $\rho \in \im u$, and this proves our claim.
\end{rk}
The following result is a partial extension of the characterization (C2) of POG curves given in the previous section to higher dimensions.
For the partial extension of the converse implication see Theorem \ref{thmPOG20} below.
\begin{thm}
\label{thmPOG} 
Let $V:f=0$ be a strictly POG hypersurface with respect to $\overline \rho=(\rho_1, \ldots, \rho_{n+1})$ with exponents
$(d_1, \ldots, d_n, d_{n+1})$. Then $d= \deg f$ is given by
$$d=d_1+ \ldots + d_n$$
and, if $\dim \Sigma=n-2$, then
$$\deg \Sigma=\sum_{j=1}^nd_j^2+\sum_{1 \leq i<j \leq n}d_id_j-\sum_{j=1}^{n+1}d_j.$$
In particular $v(\rho_n) \ne 0$ and $\deg v(\rho_n)=1$. 
\end{thm}
The equality $d=d_1+ \ldots + d_n$ for a POG hyperplane arrangement was obtained already in \cite[Proposition 4.1]{Abe}.
\proof
To prove the equality $d=d_1+ \ldots + d_n$, we use the same aproach
as in the proof of Theorem \ref{thmF0}. The equality \eqref{eqF1} has to be replaced by the exact sequence
$$0 \to S(-d_{n+1}-1) \to \oplus_{j=1}^{n+1}S(-d_j) \to D_0(f) \to 0,$$
which follows from the definition of a POG hypersurface.
This yields the following
$$\dim D_0(f)_k=\sum_{j=1}^{n+1}\binom{k+n-d_j}{n}-\binom{k+n-d_{n+1}-1}{n} .$$
Next, using the exact sequence \eqref{eqF2}, we get the following equality which replaces \eqref{eqF3}.
$$ \dim M(f)_{k+d-1}= \sum_{j=1}^{n+1}\binom{k+n-d_j}{n}-\binom{k+n-d_{n+1}-1}{n} -$$
$$  -(n+1)\binom{k+n}{n}+\binom{k+n+d-1}{n}.$$
The right hand side of this equality is shown to be a polynomial in $k$ of degree $n-1$ with leading coefficient
$$c_{n-1}=d-(d_1+\ldots + d_n).$$
By comparing to \eqref{eqD1}, we see that this leading coefficient is in fact zero. The coefficient of $k^{n-2}$ is given by
$$c_{n-2}=\frac{1}{n!}\binom{n}{2}(\sum_{j=1}^nd_j^2+d^2-2d-2d_{n+1}).$$
If $\dim \Sigma=n-2$, then we get
$$\deg \Sigma=(n-2)!c_{n-2}=\sum_{j=1}^nd_j^2+\sum_{1 \leq i<j \leq n}d_id_j-\sum_{j=1}^{n+1}d_j.$$
If $\dim \Sigma <n-2$, then
$$\sum_{j=1}^nd_j^2+\sum_{1 \leq i<j \leq n}d_id_j-\sum_{j=1}^{n+1}d_j=0,$$
and hence in this case $d_{n+1}$ is determined by the $d_j$'s with $j\leq n$.
It remains to show that $v(\rho_n) \ne 0$.
Note that $v(\rho_n)=0$ means that there is a non-trivial second order syzygy
$$\delta=a_0E+a_1\rho_1+ \ldots + a_n \rho_n=0$$
involving the Euler derivation
$$E=(x_0,x_1, \ldots, x_n)$$
and the syzygies $\rho_j \in D_0(f)$. By applying $\delta$ to $f$, we see that $a_0=0$, and hence we get a non-trivial second order syzygy
$$\delta=a_0E+a_1\rho_1+ \ldots + a_n \rho_n=0$$
involving  the syzygies $\rho_j \in D_0(f)$ for $j<n+1$, and hence not a multiple of the syzygy in \eqref{eqD2}. This contradiction proves our final claim. 
\endproof

\begin{ex}
\label{ex0} 
 Consider the surface
$$V:f=(x^3-yzw)^3+(w^3-xyz)^3=0,$$
which is a union of 3 cubic surfaces in $\PP^3$.
Then the computer algebra software SINGULAR \cite{Sing} gives $m=8$ and the
exponents $(1,4,4,7,7,8,8,8)$. In particular we have
$$d_1+d_2+d_3=1+4+4=9=\deg V,$$
but $V$ is far from being a POG surface. The Hilbert polynomial for $M(f)$ is $38t-119$, which shows that $\dim \Sigma=1$. 
The software SINGULAR \cite{Sing} gives the following first 3 generating syzygies.
$$\rho_1=(0,y,-z,0),$$
$$\rho_2=(xyz^2+3zw^3,9x^2w^2,-yz^3,3x^3z+yz^2w)$$
and
$$\rho_3=(xy^2z+3yw^3,0,-y^2z^2+9x^2w^2,3x^3y+y^2zw).$$
For both pairs $\hat \rho=(\rho_1,\rho2)$ and $\hat \rho=(\rho_1,\rho_3)$
we get $\hit(I(\hat \rho))=1<2$, hence Theorem \ref{thmT0} tells us that
$V$ is not tame with respect to these two choices for $\hat \rho$.
However, consider the last generating syzygy listed by SINGULAR, namely
$\rho_8=(A,B,C,D)$, where 
$$A=x^6yz-3x^2y^2z^2w^2+y^3z^3w^2+6x^2yzw^4+6xyzw^5-3w^8, \ B=0$$
$$C=-2x^2y^2z^3w+18x^4zw^3 \text{ and } D=3x^8-x^2y^3z^3+5x^2y^2z^2w^2+2xy^2z^2w^3-yzw^6).$$
If we take now $\hat \rho=(\rho_1, \rho_8)$, a direct computation by the computer algebra software SINGULAR shows that $\hit(I(\hat \rho)=2$, and hence $V$ is tame with respect to this $\hat \rho$ by Theorem \ref{thmT0}.
\end{ex}
The next result shows that strict POG hypersurfaces are intimately related to tame POG hypersurfaces.

\begin{prop}
\label{prop1} 
Let $V:f=0$ be a POG hypersurface with respect to $\overline \rho=(\rho_1, \ldots, \rho_{n+1})$. Then the hypersurface $V:f=0$ is strictly POG with respect to $\overline \rho$ and $\ell$ does not divide $b_n$  if and only if $V:f=0$
is tame with respect to the subset $\hat \rho=(\rho_1, \ldots, \rho_{n-1})$.
\end{prop}
\proof
Suppose first that $V$ is strictly POG and $\ell$ does not divide $b_n$ in \eqref{eqD2}. If condition \eqref{eqB2} fails,
then this is equivalent to the existence of a non trivial relation of the form
$$c_1\rho_1+ \ldots + c_{n-1}\rho_{n-1}=0.$$
Such a relation cannot be a multiple of the relation \eqref{eqD2}, since $\ell \ne 0$. Hence the condition  \eqref{eqB2} holds.
Let now 
$$\rho = \sum_{j=1}^{n+1}c_j\rho_j$$
be a syzygy not in $\im u $, and hence at least one of $c_n$ and $c_{n+1}$ is non-zero. We can assume moreover that $\ell $ does not divide $c_{n+1}$ in view of \eqref{eqD2}.
Over the field of fractions $K$ of the polynomial ring $S$, we know as above that $\rho_1, \ldots, \rho_{n-1}$ are linearly independent. If $\rho_n$ is in the $K$-vector space spanned by $\rho_1, \ldots, \rho_{n-1}$, then the $K$-vector space spanned by $\rho_1, \ldots, \rho_{n+1}$ would be of dimension $n-1$, a contradiction, since this dimension equals the rank of $D_0(f)$ and is equal to $n$.
It follows that $v(\rho_n) \ne 0$. Assume now that $v(\rho)=0$, hence
$$v(\rho)=c_nv(\rho_n)+c_{n+1}v(\rho_{n+1})=0.$$
By applying $v$ to \eqref{eqD2} we get
$$b_nv(\rho_n)+\ell v(\rho_{n+1})=0.$$
Using $v(\rho_n) \ne 0$, it follows that
$$c_n \ell= c_{n+1}b_n.$$
This is a contradiction, since $\ell$ divides neither $b_n$ nor $c_{n+1}$.
Conversely, assume that $V:f=0$ is tame. If $\ell=0$, then by taking
$\rho=\rho_n$ we see that the equality $\ker v=\im u$ does not hold.
Indeed, $\rho_n \notin \im u$, since we have a minimal set of generators.
On the other hand, the equality \eqref{eqD2} with $\ell=0$ implies that
the equality $\ker v=\im u$ fails.

\endproof

\begin{ex}
\label{ex1} 
Consider the surface
$$V:f=(x^2+y^2-zw)^2-2w^4=0.$$
Then the computer algebra software SINGULAR \cite{Sing} gives the following
syzygies
$$\rho_1=(0,w,2y,0), \ \rho_2=(-y,x,0,0), \rho_3=(w,0,2x,0)$$
and
$$\rho_4=(xz,yz,2z^2-8w^2,2x^2+2y^2-2zw).$$
There is just one generator for the second order syzygies, namely the
relation
$$x\rho_1-w\rho_2-y\rho_3=0.$$
It follows that $V:f=0$ is, for this ordering of the syzygies $\rho_j$'s, a non strictly POG surface. This shows in particular that the hypothesis strictly POG is necessary in Theorem \ref{thmPOG}. Indeed, without this hypothesis, the first exact sequence in the proof of Theorem \ref{thmPOG} might not be the minimal resolution of $D_0(f)$.
Note also that $V$ is not tame with respect to 
$\hat \rho=(\rho_1,\rho_2).$
However, if we set $\rho_1'=\rho_4$ and $\rho_4'=\rho_1$, and consider the syzygies in the order
$\rho_1', \rho_2,\rho_3, \rho_4'$, then $V:f=0$ is a strictly POG surface, and even a tame surface by Proposition \ref{prop1}, with exponents $(2,1,1,1)$. Using Theorem \ref{thmPOG}, we see that
$$\deg \Sigma=4+1+1+2+2+1-(2+1+1+1)=6.$$
This shows that we have equality in the last inequality in Theorem \ref{thmF0}, where $d_3''=1$ if we use the 3 syzygies of degree 1, even if $V$ is not free. This is due to the fact that $V$ is not tame with respect to $\hat \rho=(\rho_1,\rho_2)$, but only with respect to $\hat \rho''=(\rho_1',\rho_2)$. If we use $d_1'=2$ and $d_2'=1$, then
$d_3''=0$ and hence the inequality given by Theorem \ref{thmF0}
becomes
$$6=\deg \Sigma < 9-1 \times 2=7,$$
in accord with the fact that $V$ is not free.

\end{ex}

The proof of Proposition \ref{prop1} and Theorem \ref{thmPOG} yield the following.

\begin{cor}
\label{cor1} 
If $V:f=0$ is a non strictly POG with respect to $\overline \rho=(\rho_1, \ldots, \rho_{n+1})$, 
 then one of the following two cases occurs.
\begin{enumerate}

\item If  $b_n\ne 0$, then $B(f,\hat \rho)$ is a principal ideal generated by the polynomial $v(\rho_{n+1})$.

\item  If  $b_n= 0$, then $B(f,\hat \rho)=0$.

\end{enumerate}

If $V:f=0$ is a strictly POG with respect to $\overline \rho=(\rho_1, \ldots, \rho_{n+1})$, then one of the following two cases occurs.
\begin{enumerate}

\item [(i)]
 The Bourbaki ideal $B(f,\hat \rho)$ defines a codimension 2 complete intersection of multi-degree $(1,d_{n+1}-d_n+1)$ if  $V$ is tame with respect to the subset $\hat \rho$.

\item [(ii)]
 The Bourbaki ideal $B(f,\hat \rho)$ defines a hyperplane if  $V$ is not tame with respect to the subset $\hat \rho$.

\end{enumerate}
In particular, a POG hypersurface as above is tame if and only if the corresponding Bourbaki ideal $B(f,\hat \rho)$ defines a codimension 2 complete intersection of multi-degree $$(1,d_{n+1}-d_n+1).$$

\end{cor}
The last claim in Corollary \ref{cor1} may be regarded as a refinement of the last claim in Theorem \ref{thmT1} saying that $\gr (B(f,\hat \rho))=2$.
\proof
Consider first the case of a non strictly POG such that $b_n \ne 0$.
This implies $v(\rho_n)=0$.
On the other hand, note that $v(\rho_{n+1}) \ne 0$, since otherwise we get a second order syzygy
$$\delta=a_0E+a_1\rho_1+ \ldots + a_{n-1} \rho_{n-1}+a_{n+1} \rho_{n+1}  =0$$
as in the proof of Theorem \ref{thmPOG} above. As above, we have $a_0=0$ and this leads to a contradiction with the assumption $b_n \ne 0$.
Consider now the case of a non strictly POG such that $b_n = 0$. Then the syzygies in $\hat \rho$ are linearly dependent over the fields of fractions $K$, and hence the map $v$ is the trivial map in this case.

Assume now that  is a strictly POG as in \eqref{eqD2}. Then we know that $\ell_1=v(\rho_n)$ is a non-zero linear form by Theorem \ref{thmPOG}. If $\ell$ does not divide $b_n$, then $\ell$ divides $\ell_1$, and $v(\rho_{n+1})$ is just $b_n$ up to a nonzero factor. Hence $B(f,\hat \rho)=(\ell,b_n)$ in this case.

If $\ell$ divides $b_n$, then $\ell_1$ divides $v(\rho_{n+1})$ and hence
$B(f,\hat \rho)=(\ell_1)$ in this case.

\endproof

\begin{ex}
\label{ex2} 
Consider the surface
$$V:f=(x+y+z)(x+2y)wz+xyzw=0$$
which is a union of two planes and a quadratic surface in $\PP^3$.
Then the computer algebra software SINGULAR \cite{Sing} shows that
$m=4$ and the exponents are $(1,1,2,2)$. The Bourbaki ideal corresponding to the choice $\hat \rho=(\rho_1,\rho_2)$, the two syzygies of degree 1, is given by
$$B(f,\hat \rho)=(y,2x+z),$$
and the minimal resolution of the Jacobian ideal is
$$0 \to S(-6) \to S(-4)^2 \oplus  S(-5)^2 \to S(-3)^4 --> S.$$
It follows that $V$ is a tame, nearly free surface.

\end{ex}

\begin{thm}
\label{thmPOG20} 
Let $V:f=0$ be a tame hypersurface with respect to
$\hat \rho=(\rho_1, \ldots, \rho_{n-1})$, which is an $(n+1)$-syzygy hypersurface with exponents
$(d_1, \ldots, d_{n+1})$ such that $d= \deg f$ is given by
$$d=d_1+ \ldots + d_n.$$
Then $V$ is a strictly POG hypersurface with respect to $(\rho_1, \ldots, \rho_{n+1})$. 
\end{thm}
\proof
It follows from Remark \ref{rkT0} that the scheme $Z$ defined by the Bourbaki ideal $B(f, \hat \rho)$ has codimension 2.
Exactly as in the final step of the proof of Theorem \ref{thmPOG}, we see that $\ell_1=v(\rho_n)$ is a non-zero linear form. Since $V$ is 
an $(n+1)$-syzygy hypersurface, it follows that the second order syzygies of $D_0(f)$ are generated by a unique relation
\begin{equation}
\label{eqsyz2}
a_1\rho_1+ \ldots + a_{n-1}\rho_{n-1}+a_{n}\rho_{n}+a_{n+1}\rho_{n+1}=0,
\end{equation}
where both $a_n$ and $a_{n+1}$ are non-zero polynomials of degree $>0$. This relation implies
$$a_{n}v(\rho_{n})+a_{n+1}v(\rho_{n+1})=0.$$
Since $B(f, \hat \rho)=(v(\rho_{n}),v(\rho_{n+1}))$ has codimension 2, if $\deg a_{n+1}>1$, then $a_{n+1}=\ell_1b$ and $a_n=-bv(\rho_{n+1})$, with $\deg b >0$. This implies that
$$a_1\rho_1+ \ldots + a_{n-1}\rho_{n-1}-b(v(\rho_{n+1})\rho_{n}+v(\rho_{n})\rho_{n+1})=0,$$
and hence 
$$v(\rho_{n+1})\rho_{n}+v(\rho_{n})\rho_{n+1}\in \ker v=\im u.$$
This contradicts the fact that \eqref{eqsyz2} generates the second order syzygies of $D_0(f)$, and hence proves our claim.
\endproof

\section{Some questions related to Jacobian syzygies of hypersurfaces} 

\begin{question}
\label{q1}
Find an upper bound for $m$, the minimal number of generators for the graded $S$-module $D_0(f)$ in terms of $d$ and $n$. 
When $V:f=0$ is a generic hyperplane arrangement in $\PP^n$ then it follows from \cite[Theorem 4.5.3 and Corollary 4.5.4]{RT} that one has
$$m=\binom{d-1}{n-1}.$$
Does the bound
\begin{equation}
\label{b1}
m\leq \binom{d-1}{n-1}
\end{equation}
hold for any arrangement of $d$ hyperplanes in $\PP^n$ ?
The answer is positive for $n=2$, see \eqref{e2}.
When $V$ is a degree $d$ hypersurface in $\PP^n$ having only isolated singularities, the following bound
\begin{equation}
\label{b2}
m \leq \binom{n+1}{2}+ \deg(\Sigma)
\end{equation}
seems to  work. Note that we have equality for a smooth hypersurface $V:f=0$ in $\PP^n$. Moreover, when $V$ has only isolated singularities, one has
$$\deg (\Sigma )= \tau(V),$$
the total Tjurina number of $V$,
and the above upper bounds seems to hold,  compare to \cite[Corollary 11]{CD}. Note that the homology group $H_1(K)$ which appears in that paper is, up to a shift, the quotient of $D_0(f)$ by the $S$-submodule
generated by all Koszul type derivations of $f$, and there are clearly
$ \binom{n+1}{2}$ such derivations, all of them of degree $(d-1)$.

It remains of course the problem to bound $\deg(\Sigma)$ in terms of
$n$ and $d$. When $V$ has at most isolated singularities one has
$$\deg (\Sigma )= \tau(V) \leq (d-1)^n.$$
On the other hand, when $V$ is a hyperplane arrangement in $\PP^n$,
and $V'$ is the intersection of $V$ with a generic 2-plane $E$ in $\PP^n$, it follows \cite[Remark 2.3]{WaYo} that one has
$$\deg (\Sigma )= \tau(V') \leq (d-1)^2.$$
It follows that the two upper bounds \eqref{b1} and \eqref{b2} are quite different in nature. 
\end{question} 

\begin{ex}
\label{ex110} 
Consider the surface
$$V:f=x^4 +y^4+z^4 - xyzw=0.$$
A direct computation using SINGULAR shows that $m=12$ and the exponents are $( 3, 3, 3, 3, 3, 3, 4, 4, 4, 4, 4, 4 )$.
Note that $V$ has isolated singularities and
$$\tau(V)= \deg (\Sigma)=  10,$$
and hence the inequality
$$12 \leq 6+10$$
suggested in Question \ref{q1} holds.

\end{ex}

All examples we have computed so far suggest that the following question has a positive answer
\begin{question}
\label{q2}
Find an upper bound for the maximal degree $d_m$ of a Jacobian syzygy in terms of $d$.
Does the bound
$$d_m \leq n(d-2),$$
compare to \eqref{e3}, 
works for any $n$ ? We know that this bound works when $\dim \Sigma=0$, see  \cite[Corollary 11]{CD}.
\end{question} 

Concerning the tame hypersurfaces, the following two key questions are 
related.
\begin{question}
\label{q2.5}
 Does there exist hypersurfaces $V:f=0$ which are not tame for any choice of a subset
$\hat \rho$ of $n-1$ generators in a minimal set of generators for the graded $S$-module $D_0(f)$ ?

\end{question} 

Theorem \ref{thmT1} says that any hypersurface $V:f=0$ is tame with respect to $\hat \rho'=(\rho_1', \ldots, \rho_{n-1}')$, where 
$\deg \rho_j'=d_m$, for any $j=1,\ldots,n-1$. Such a subset 
$\hat \rho'$ is not a subset of $n-1$ generators in a minimal set of generators for the graded $S$-module $D_0(f)$. 
If the answer to the Question \ref{q2.5} is negative, then for any hypersurface $V$ we can define a new invariant $b(V)$, namely the minimal degree of the subset
$$\deg (\hat \rho)=d_1+\ldots + d_{n-1},$$
where $\hat \rho=(\rho_1, \ldots, \rho_{n-1})$, $d_j= \deg\rho_j$
and $\hat \rho$ runs through the set of subsets of generators such that
$V$ is tame with respect to $\hat \rho$.

We have seen in Proposition \ref{prop1} and in Example \ref{ex1} that being a strictly POG hypersurface depends on the choice of the first $(n-1)$ syzygies in $\hat \rho$. A similar example is the following.

\begin{ex}
\label{ex10} 
 Consider the 3-fold
$$V:f=(x^2+y^2+z^2)(u^2+v^2)=0,$$
which is a union of two planes and a quadratic cone in $\PP^4$.
Then the computer algebra software SINGULAR \cite{Sing} gives $m=5$ and the
exponents $(1,1,1,1,1)$, namely one can take
$$\rho_1=(0,0,0,-v,u), \  \rho_2=(0,-z,y,0,0), \ \rho_3=(x,y,z,u,v),$$
$$\rho_4=(-y,x,0,0,0) \text{ and } \rho_5=(-z,0,x,0,0).$$
The minimal resolution of the Jacobian ideal is
$$0 \to S(-5)  \to  S(-4)^5 \to S(-3)^5  \to  S,$$
hence $V$ is a nearly free 3-fold. On the other hand, we have the following second order syzygy
$$x\rho_2+z\rho_4-y\rho_5=0.$$
It follows that for the choice $\hat \rho=(\rho_2,\rho_4, \rho_5)$, the 3-fold $V$ is not a strictly POG hypersurface and the corresponding Bourbaki ideal $B(f,\hat \rho)$ is trivial. But for the choice
$\hat \rho '=(\rho_1,\rho_2,\rho_3)$, the 3-fold $V$ is  a strictly POG hypersurface.
\end{ex}

\begin{question}
\label{q3}
1) Does there exist a POG hypersurface  $V:f=0$ with respect to $\overline \rho=(\rho_1, \ldots, \rho_{n+1})$ which is not strictly POG for any choice of a subset
$\hat \rho$ of $n-1$ generators in a minimal set of generators for the graded $S$-module $D_0(f)$

2) Does there exist a POG hypersurface  with respect to $\overline \rho=(\rho_1, \ldots, \rho_{n+1})$ which is not tame with respect for any choice of a subset
$\hat \rho$ of $n-1$ generators in a minimal set of generators for the graded $S$-module $D_0(f)$

3) Same questions as 1) and 2), but restricted to the class of nearly free hypersurfaces.

\end{question} 

In relation with the example of a tame hypersuface $V$ with $\dim \Sigma <n-2$ given in Remark \ref{rkF0}, one may ask the following.
A negative answer to question 1) below would simplify the hypotheses of Theorem \ref{thmPOG} as well.
\begin{question}
\label{q4}
1) Does there exist a  POG hypersurface such that $\dim \Sigma <n-2$ ?

2) Does there exist a  tame POG hypersurface such that $\dim \Sigma <n-2$ ?

\end{question}

\end{document}